\documentclass{article}

\begin{document}

\title{Stability of algebraic varieties and K\"ahler geometry}
\author{Simon Donaldson}
\date{\today}
\maketitle



\newcommand{\cL}{{\cal L}}
\newcommand{\cF}{{\cal F}}
\newcommand{\ob}{\overline{b}}
\newcommand{\oz}{\overline{z}}
\newcommand{\dbd}{\overline{\partial}\partial}
\newcommand{\bR}{{\bf R}}
\newcommand{\bC}{{\bf C}}
\newcommand{\bZ}{{\bf Z}}
\newcommand{\bP}{{\bf P}}
\newcommand{\dbar}{\overline{\partial}}
\newtheorem{defn}{Definition}

\section{ Background}

\subsection{K\"ahler metrics}

This write-up follows quite closely the author's three lectures at the AMS Summer Research Institute. The context for the lectures was the search for  \lq\lq optimal'' Riemannian metrics on compact complex manifolds and in particular on complex projective manifolds. This is a huge subject and we do not attempt to cover everything; we focus on situations in which these differential geometric questions are related to Geometric Invariant Theory and notions of \lq\lq stability''. 
We begin by recalling some standard background.   
Let $X$ be an $n$-dimensional complex manifold. A hermitian metric on $X$ is  a Riemannian
metric  which is algebraically compatible with the complex structure,
in the sense that it is the real part of a Hermitian form. The imaginary
part of the Hermitian form is skew-symmetric, an exterior $2$-form $\omega$.
 The metric is {\it K\"ahler} if $\omega$ is closed. There are several equivalent definitions.  One is that the holonomy of the Riemannian metric lies
in the unitary group---i.e. parallel transport of tangent vectors commutes with multiplication
by $I$. Another is that at each point there is a complex co-ordinate system in which the metric agrees with the flat metric to first order. In any case, this K\"ahler condition is the natural notion of compatibility between the complex and metric structures and we will always consider K\"ahler metrics here. (K\"ahler geometry can be seen as the intersection of the three branches:  complex, Riemannian
and symplectic geometries.)

In local co-ordinates $z_{a}$  the Hermitian structure is given by a matrix
$g_{a\ob}$ and $\omega= i\sum g_{a\ob} dz_{a} d\oz_{\ob}$. The $ \dbd$-Lemma
 states that the metric can be given, locally, by a potential $\phi$

$$    g_{a\ob} = \frac{\partial^{2}\phi}{\partial z_{a}\partial \oz_{\ob}},
$$  
or, in  different notation,  $\omega= i \dbd \phi$.

\

Suppose that $X$ is compact and $\omega_{0}$ is a K\"ahler metric on $X$.
It defines a non-zero cohomology class $[\omega_{0}]\in H^{2}(X;\bR)$. The
global version of the $\dbd$-Lemma says that any other K\"ahler metric $\omega$
in the same cohomology class can be represented by a potential
$\omega= \omega_{0} + i \dbd \phi$. This has a geometric interpretation in the case when $2 \pi [\omega_{0}]$ is the first Chern class of a holomorphic line bundle $L\rightarrow X$ (so $X$ is a projective algebraic variety, by the Kodaira embedding theorem). Recall that if $h$ is a Hermitian metric on $L$ there is a unique unitary connection on $L$ compatible with the complex structure.  The form $ -2\pi i \omega_{0}$ can then be realised as the curvature of a reference metric $h_{0}$ (i.e. the curvature of the compatible connection), and $\omega$ similarly corresponds to a metric $e^{\phi} h_{0}$.

Recall that a Riemannian manifold has a curvature tensor ${\rm Riem}$. This
has a contraction to the Ricci tensor,  which is closely related to the induced volume form. One geometric characterisation
of the Ricci tensor is in terms of the Riemannian volume form in local geodesic
co-ordinates $x_{i}$

  \begin{equation} {\rm Vol} =  (1 - \frac{1}{6} \sum_{ij} R_{ij}x_{i}x_{j}+ O(x^{3})) {\rm
Vol}_{{\rm Euc}}. \end{equation} 

A further contraction of the Ricci tensor gives the {\it scalar curvature}
$S$.

On a K\"ahler manifold the Ricci tensor is Hermitian,  the corresponding
$2$-form $\rho$ is closed and represents the class $2\pi c_{1}(X)$. In local
complex co-ordinates, a formula for the Ricci form is
$$       \rho= i \dbd {\rm \log(Vol)}, $$
where ${\rm Vol}= \det ( h_{a\ob})$, the volume in local co-ordinates.

We consider four notions of \lq\lq optimal'' metrics.
\begin{itemize} \item If $c_{1}(X)$ is a multiple of the K\"ahler class,
$c_{1}= \lambda [\omega]$, then we can look for a {\it K\"ahler-Einstein}
metric: $ \rho= \lambda \omega$. 
\item In any K\"ahler class we may seek a constant scalar curvature K\"ahler
(CSCK) metric. Clearly K\"ahler-Einstein metrics are CSCK. 
\item In any K\"ahler class we may consider the critical points of the
{\em Calabi functional} $ \Vert {\rm Riem} \Vert^{2}_{L^{2}} $.  These are
called
{\it extremal metrics}. The Euler-Lagrange equation is that the
gradient of the  scalar curvature is a holomorphic vector field. 
CSCK metrics are extremal (and on a manifold with no holomorphic vector fields
the converse is true).
\item In the case when $c_{1}(X)=\lambda [\omega]$, the  {\it K\"ahler-Ricci
soliton} equation is  that $\rho-\lambda \omega$
is the Lie derivative of $\omega$ along a holomorphic vector field. Clearly K\"ahler-Einstein metrics are solitons. 
\end{itemize}

Extremal metrics and K\"ahler-Ricci solitons arise naturally when considering associated evolution equations, for 1-parameter families of metrics. The {\it Calabi flow} is the equation 
$$   \frac{\partial \omega}{\partial t} = i \dbd S(\omega), $$
while the {\it K\"ahler-Ricci flow} is 
$$  \frac{\partial \omega}{\partial t} = (\lambda \omega-\rho), $$
(in the case when $c_{1}(X)= \lambda[ \omega]$). Extremal metrics are fixed points of the Calabi flow modulo holomorphic automorphisms and likewise for solitons and the K\"ahler-Ricci flow. 
  
  \
  
  \

The basic question we want to address then is the existence of these metrics, in a given K\"ahler class on a  given complex manifold $X$.  There is ample motivation for studying this problem. Of course a paradigm is the case of complex dimension 1,  where a compact Riemann surface
admits a metric  of constant Gauss curvature. This statement can be viewed as a formulation of
the uniformisation theorem for compact Riemann surfaces. In higher dimensions there is a huge literature and renowned general theorems; notably from the 1970's the
existence of K\"ahler-Einstein metrics in the negative case ($\lambda<0$)
\cite{kn:Aub}, \cite{kn:Y}
 and the Calabi-Yau case ($\lambda=0$) \cite{kn:Y} . The existence theorems have important consequences for {\it  moduli problems}, which we return to in 4.2 below.  The over-riding motivation is perhaps that these existence problems lead to fascinating interactions between Riemannian
geometry,  complex (algebraic) geometry, PDE and analysis.

\subsection{Geometric Invariant theory, Kempf-Ness etc.}

Going back to the 1980's, it has been realised that there is a general package
of ideas relating complex moduli problems, in the setting of  Mumford's Geometric Invariant Theory,  to \lq\lq metric geometry''. In this subsection we  review  this and explain why it gives a conceptual framework for the existence questions in K\"ahler geometry introduced above. 

Let $K^{c}$ be reductive complex Lie group, the complexification of a compact
group $K$ (for example $K=SU(m), K^{c}= SL(m,\bC)$). Let $V$ be a complex representation
of $K^{c}$. Then we have a \lq\lq complex moduli problem''--- to study the set of orbits of $K^{c}$
in
$\bP(V)$. The fundamental difficulty is that the set of all orbits will not have
a good structure---it will not be Hausdorff in the induced topology.  The algebraic approach is to consider the graded  ring $R$ of $K^{c}$-invariant polynomials on $V$. A famous theorem of Hilbert states that this is finitely generated,
so ${\rm Proj}(R)$ is a  projective variety. The points  of ${\rm Proj}(R)$
correspond approximately, but not exactly to the $K^{c}$ orbits in $\bP(V)$ and we write
$$   {\rm Proj}(R)= \bP(V)//K^{c}. $$

 We say that a nonzero $v\in V$ is \begin{itemize}
 \item {\it polystable} if its $K^{c}$-orbit is closed; \item {\it stable} if it is polystable and its stabiliser in $K^{c}$ is finite;
\item  {\it semi-stable} if $0$ does not lie in the closure
of its $K^{c}$-orbit. \end{itemize} These notions pass down to the projectivisation and  basic result of Mumford's theory \cite{kn:Mum} is that 
the points of $\bP(V)//K^{c}$ correspond to the polystable orbits.

Now introduce \lq\lq metric geometry'' by fixing a $K$-invariant Hermitian
metric on $V$. The beautiful and elementary observation of Kempf and Ness \cite{kn:KN} is that a point $v\in V\setminus \{0\}$ is polystable  if and only if its $K^{c}$-orbit contains a point
which minimises the norm within the orbit, and this minimiser is unique up to the action of $K$. Semistability is characterised by the property that the norm functional has a strictly positive lower bound.

 Given a non-zero point $v\in V$ the function ${\cal F}_{v}(g)= \log \vert
g(v)\vert$ on $K^{c}$ is preserved by left multiplication by $K$ so descends
to a function
$$  {\cal F}_{v}: K^{c}/K\rightarrow \bR. $$

The manifold $K^{c}/K$ is a {\it symmetric space of negative curvature}.
The stability or semi-stability of $v$ depends on whether ${\cal F}_{v}$ attains
a minimum or is bounded below. An important general property is that ${\cal
F}_{v}$ is {\it convex} along geodesics in $K^{c}/K$. This gives the uniqueness
property and also shows that any critical point of the norm functional on an orbit is a minimum.   There is a map
$$  \mu : V\rightarrow {\rm Lie}(K)^{*}$$
 (the dual of the Lie algebra) defined by the first order variation of the norm. That is, for $v$ in $V$ and $\xi$ in ${\rm Lie}(K)$ 
$$  \frac{d}{dt}\vert_{t=0}\left( \vert  \exp(i \xi) v \vert\right)= \langle \mu(v), \xi\rangle. $$
The map $\mu$ his homogeneous of degree $2$ so $\mu_{\bP}(v)= \vert v \vert^{-2} \mu(v)$ factors through $\bP(V)$ Then we have
   $$  \bP(V)//K^{c}= \mu_{\bP}^{-1}(0)/K. $$
 \

 The {\it Hilbert-Mumford numerical criterion} gives a practical way to detect stability. For simplicity suppose that
the
$K^{c}$-orbit of $v$ is free. Let $\lambda: \bC^{*}\rightarrow K^{c}$ be
a non-trivial 1-parameter subgroup. Then $\lambda(t)v$ has a Laurent expansion
about $t=0$. Let $-w(\lambda,v)\in \bZ $ be the order of the leading term
in this expansion. Then $v$ is stable if and only if $w(\lambda,v)>0$ for
all  $\lambda$ and semi-stable if weak inequality holds.

We will illustrate these ideas with two simple examples. 

\

{\bf Example 1} 
Take $K=SU(2), K^{c}=SL(2,\bC)$  and let $V=U^{\otimes d}$ be the d-fold tensor product of the standard 2-dimensional representation $U=\bC^{2}$. Restrict attention to the subset $W\subset V\setminus \{0\}$ of non-zero decomposable tensors. Thus the projectivisation of $W$ can be identified with the d-fold product $(\bC\bP^{1})^{d}$. One finds  that for $u_{1},\dots u_{d}$ in $\bP^{1}$ 
$$  \mu_{\bP}(u_{1}, \dots ,u_{d})= \sum_{i=1}^{d} m(u_{i}), $$
where $m:\bC\bP^{1}\rightarrow {\rm Lie}\ K= \bR^{3}$ is the standard embedding of $\bC\bP^{1}$ as the unit 2-sphere. Thus $\mu_{\bP}^{-1}(0)$ is identified with the configurations of $d$ ordered points on the sphere which have centre of mass at $0$. The Hilbert-Mumford criterion shows that a configuration $(u_{1},\dots, u_{d})$ is polystable if either no point occurs with multiplicity  $\geq d/2$ or if $d$ is even and there are distinct points each of multiplicity $d/2$. For example when $d=4$ one has $\bP(W)//K^{c}= \bC\bP^{1}$ with the equivalence defined by the cross-ratio.  The stability discussion corresponds to the fact that one can define the cross-ratio of a quadruple with  for configurations with multiplicity  2,  but not 3. 

\

 {\bf Example 2}
 Take $K=U(n), K^{c}=GL(n,\bC)$ and let $V= {\rm End}(\bC^{n})$ be the adjoint representation. Using the standard identification of the Lie algebra and its dual, one finds that
$  \mu(A)= i[A,A^{*}]$,  where $A^{*}$ denotes the usual adjoint. 
So $\mu^{-1}(0)$ consists of the \lq\lq normal'' matrices, which commute with their adjoints. On the complex side, the polystable points are the diagonalisable matrices and the identification $\bP(V)//K^{c}= \mu_{\bP}^{-1}(0)/K$ boils down to the familiar facts that  normal matrices are  diagonalisable and conversely a diagonalisable matrix is conjugate to a normal one (just map the eigenvectors to an orthonormal basis).

 \

 The discussion above is in finite dimensions but 
 many examples (often going under the general heading of the \lq\lq Kobayashi-Hitchin
correspondence'') are understood of {\em differential geometric} problems
which can be viewed as an  {\em infinite-dimensional versions} of this set-up. A key step in this is a symplectic geometry point of view. In the setting above, we focus on the action of $K$ on the projectivisation $\bP(V)$, preserving the standard Fubini-Study metric and the corresponding symplectic form $\omega_{FS}$. Then the map $\mu_{\bP}:\bP(V)\rightarrow {\rm Lie}(K)^{*}$ appears as the \lq\lq moment map'' for the action. That is, for each $\xi\in {\rm Lie}(K)$ the component $\langle \mu_{\bP}(\ \ ), \xi\rangle$ is a Hamiltonian function for the 1-parameter subgroup generated by $\xi$. One can then recover the norm function $F$ on $K^{c}/K$ from $\mu$ by integration. From this point of view  the essential ingredients are:
\begin{itemize} \item A complex  manifold ${\cal P}$ with a K\"ahler form $\Omega$.
\item A group ${\cal G}$ with an isometric holomorphic action on ${\cal P}$.
\item A moment map $\mu: {\cal P}\rightarrow {\rm Lie}({\cal G})^{*}. $
\item A complexified group ${\cal G}^{c}\supset {\cal G}$ and an extension of the action to a holomorphic action of ${\cal G}^{c}$ on ${\cal P}$. 
\end{itemize}
Given such as set-up we can ask for a useful definition of \lq\lq polystable points'' in ${\cal P}$ and try to establish a correspondence between the resulting complex quotient ${\cal P}//{\cal G}^{c}$ and the symplectic quotient $\mu^{-1}(0)/{\cal G}$. If ${\cal P}$ is compact (in particular finite-dimensional) and the cohomology class $2\pi [\Omega]$ is integral we are almost back in the previous case, because we can embed ${\cal P}$ in projective space by the Kodaira embedding theorem, but this formulation allows us to consider infinite-dimensional situations of a fundamentally different nature.

 The first example of this kind is due to Atiyah and Bott \cite{kn:AB}.  Let $E\rightarrow \Sigma$ be a (trivial)
$C^{\infty}$-complex vector bundle over a compact Riemann surface. Let ${\cal
A}$ be the set of $\dbar$-operators on $E$: an infinite dimensional complex affine space. There is a group ${\cal G}^{c}$
of bundle automorphisms of $E$ and the orbits of ${\cal G}^{c}$ in ${\cal
A}$ correspond to isomorphism classes of {\it holomorphic} vector bundles.
 Now fix a Hermitian metric on $E$. We get a subgroup ${\cal G}\subset {\cal
G}^{c}$ of unitary automorphisms. Points in ${\cal A}$ can be viewed as unitary
connections on $E$. Tangent vectors to ${\cal A}$ are $1$-forms with values in the bundle ${\rm End} E$ and there is a symplectic form given by
$$  \Omega(a,b)= \int_{\Sigma} {\rm Tr} (a\wedge b), $$
 which is a K\"ahler form with respect to the complex structure. The moment map is given by the curvature of the connection so $\mu^{-1}(0)/{\cal G}$ consists of the equivalences classes of flat unitary connections. 

The case when the bundle $E$ has rank one is a piece of classical theory. In that case one has an identification ${\cal A}/{\cal G}^{c}= \mu^{-1}(0)/{\cal G}$ which is the statement that any holomorphic line bundle of degree $0$ has a unique compatible flat unitary structure. This amounts to solving a linear PDE and the notion of stability does not enter. The statement can be viewed as the identification of the Jacobian with the torus
$    H^{1}(\Sigma; \bR)/ H^{1}(\Sigma, \bZ)$.

For bundles of higher rank the situation is more interesting. On the one hand, Mumford gave a definition of a stable bundle $E\rightarrow \Sigma$. He was lead to this by embedding the moduli problem for bundles in a linear set-up of the kind we considered above. The numerical criterion which emerged is that a bundle $E$ (of degree $0$, for simplicity) is defined to be stable if  all proper sub-bundles $E'\subset E$ have strictly negative degree. Then a polystable bundle (of degree $0$) is a direct sum of stable bundles (each of degree $0$). Taking this as the definition of a polystable point in ${\cal A}$, the identification ${\cal A}//{\cal G}^{c}= \mu^{-1}(0)/{\cal G}$ is true and amounts to the 1965 theorem of Narasimhan and Seshadri \cite{kn:NS}: a bundle admits a flat unitary connection if and only if it is polystable.  

There are many other examples in this vein, and a large literature. These include the Hermitian Yang-Mills equation for bundles over higher dimensional base manifolds,  \lq\lq pairs'' consisting of a bundle with a holomorphic section and \lq\lq parabolic structures'' along a divisor in the base. In each case a necessary and sufficient condition for the existence of some differential geometric structure is given by an algebro-geometric stability condition, involving a numerical criterion. A point to emphasise is that, while there is a   commmon  conceptual framework to the results, this  framework does not in itself provide {\it proofs}. In the prototype finite-dimensional case of Kempf and Ness the proof is elementary because we can minimise the norm functional over a compact set, but this compactness is lost in infinite dimensions. The statements bear on the solubility of nonlinear PDE and the proofs require detailed analysis in each case.

We now return to our main theme of existence problems for K\"ahler metrics. These can, to some extent, be
fitted into  the conceptual framework sketched above \cite{kn:SKD-1}\cite{kn:SKD0}.
Let $(M,\omega)$ be a compact symplectic manifold and let ${\cal G}$ be the group
of \lq\lq exact'' symplectomorphisms. The Lie algebra of ${\cal G}$ is the space of
Hamiltonian functions on $M$, modulo constants. (If $H^{1}(M;\bR)=0$ then the exactness condition is vacuous.) 

Let ${\cal J}$ be the space of  almost-complex structures compatible
with $\omega$. Then certainly ${\cal G}$ acts on ${\cal J}$. Tangent vectors to ${\cal J}$ at a given almost-complex structure $J_{0}$ can be identified with sections of the complex vector bundle $s^{2}(TM)$ (using the almost complex structure $J_{0}$). As such, the tangent space acquires a complex structure and this makes ${\cal J}$, at least formally,  into an infinite-dimensional complex manifold. Similarly, the natural $L^{2}$ metric on these tensor fields gives ${\cal J}$ a K\"ahler structure. For our purposes, we are interested in the ${\cal G}$-invariant subset ${\cal J}_{\rm int}\subset {\cal J}$ of {\it integrable} almost complex structures. This is, formally,  a complex analytic subvariety of ${\cal J}$ and the possible singularities are not relevant to this discussion. So in sum we get an action of ${\cal G}$ on an infinite dimensional K\"ahler space $({\cal J}_{\rm int}, \Omega)$ and we can ask for a moment map for this action. Calculation reveals that this moment map is simply the scalar curvature $\mu(J)= S(\omega, J)$ of the metric induced by $\omega$ and $J$ \cite{kn:Fujiki}. More precisely,  we use the pairing between functions modulo constants and $S$ given by
$$     f\mapsto \int_{M} f (S-\hat{S}) \frac{\omega^{n}}{n!}, $$
where $\hat{S}$ is the average value of $S$  (which is a topological quantity determined  by $[\omega]$ and $c_{1}(M)$). So the zeros of the moment map are K\"ahler metrics of constant scalar curvature.

 We now encounter the fundamental difficulty that there is  no complexification ${\cal G}^{c}$ of the symplectomorphism group ${\cal G}$. However we have an equivalence relation on ${\cal J}_{{\rm int}}$
  given by  $J_{1}\sim J_{2}$
if $(M,J_{1}), (M,J_{2})$ are isomorphic complex manifolds. 
The equivalence classes can be interpreted as the orbits of the (nonexistent)
group ${\cal G}^{c}$. To explain this briefly, consider for simplicity the case when $(M,J,\omega)$ has no automorphisms, so the ${\cal G}$-orbit of $J$ is free. Let ${\cal O}\subset {\cal J}$ be the equivalence class of $J$, as defined above. Then one finds that the tangent bundle of ${\cal O}$ can be trivialised as ${\cal O}\times H$ where $H$ is the space of complex valued functions on $M$ modulo constants. Thus $H$ has a Lie algebra structure, regarded as the complexification of ${\rm Lie}({\cal G})$,  with the  bracket given by the complexified  Poisson bracket $\{\ ,\ \}$. Under the trivialisation, any element $\eta\in H$ defines a vector field $V_{\eta}$ on ${\cal O}$ and one has
   $$   [V_{\eta_{1}}, V_{\eta_{2}}] = V_{\{\eta_{1}, \eta_{2}\}}. $$
   If ${\cal O}$ were a finite-dimensional manifold this would lead to a Lie group structure, with the $V_{\eta}$ the left-invariant vector fields, but the relevant integration results fail in infinite dimensions. But we could say that ${\cal O}$ behaves like a group orbit for the purposes of any infinitesimal calculation. In essence, the Hamiltonian functions which make up the Lie algebra of ${\cal G}$ complexify to  K\"ahler potentials. (A detailed development of this point of view is given in \cite{kn:SKD0}.)

\

By a theorem of Moser, the symplectic structures defined by different K\"ahler metrics on a complex manifold  in a given cohomology class are all equivalent.  Thus,  given the interpretation above of the equivalence classes in ${\cal J}_{{\rm int}}$ as \lq\lq complex  orbits'', the question of the existence of a constant scalar curvature metric in a given K\"ahler class becomes precisely our \lq\lq standard question'' of finding a zero of the moment map in  a complex orbit, which we expect to be related to a notion of stability.

Taking a slightly different point of view, fix a complex manifold $X$  and let ${\cal H}$ be the space of K\"ahler metrics in a given cohomology class. In the formal picture above this is viewed as ${\cal G}^{c}/{\cal G}$. 
and indeed ${\cal H}$ does have the structure of an infinite dimensional
symmetric space of negative curvature, defined by the {\it Mabuchi metric} \cite{kn:Mabuchi2}. This is exactly analogous to the standard finite-dimensional case, when $K^{c}/K$ has such a structure. 
For an infinitesimal variation $\delta \omega = i\dbd (\delta \phi)$ (normalised
so that
  the integral of $\delta \phi$ is zero) we set:
   $$  \Vert \delta \omega \Vert^{2}_{\omega}= \int_{X} (\delta \phi)^{2}
\frac{\omega^{n}}{n!}. $$
Thus $\delta \omega$ is viewed as a tangent vector of ${\cal H}$ at $\omega$ and this formula defines a Riemannian metric on ${\cal H}$. 
The functional on ${\cal H}$ analogous to the Kemp-Ness function $\log \vert v\vert$
is the {\it Mabuchi functional} \cite{kn:Mabuchi}. It can be defined by its infinitesimal variation

$$ \delta {\cal F}= \int_{X} \delta \phi (S-\hat{S}) \frac{\omega^{n}}{n!} , $$
where $S$ is the scalar curvature and ${\hat S}$ is the average of $S$ over
$X$. 
So a critical point of ${\cal F}$ is exactly a CSCK metric.

 While we will not go into detail here, similar things hold for the other equations (extremal, Kahler-Einstein, Kahler-Ricci soliton).  In the general context of a K\"ahler manifold ${\cal P}$ with moment map $\mu: A\rightarrow {\rm Lie}(K)^{*}$ we suppose that the Lie algebra $K$ has an invariant bilinear form, so we can identify it with its dual. So $\mu(x)$ is now thought of as an element of the Lie algebra and the derivative of the action gives a vector field $ V_{\mu(x)}$ on ${\cal P}$. Then we can consider points $x\in  {\cal P}$ such that $V_{\mu(x)}(x)=0$.  These are the analogues, in the general context, of extremal metrics. For example in the case of configurations of $d$ points on the sphere the condition becomes that the configuration has the form of a point $p$ taken with some multiplicity $r$ and the antipodal point $-p$ taken with multiplicity $d-r$. There is a similar  interpretation of the Calabi flow as the gradient flow of the function $\vert \mu\vert^{2}$, as studied in the finite dimensional situation by Kirwan \cite{kn:FCK}.  K\"ahler-Einstein metrics can be regarded as special cases of constant scalar curvature ones. (An integral identity shows that when $c_{1}(X)=\lambda[\omega]$ any CSCK metric is K\"ahler-Einstein.) They can also be fitted into the general framework by using a  different K\"ahler structure on ${\cal J}_{{\rm int}}$ \cite{kn:SKD4} with another important functional, the Ding functional, taking the place of the Mabuchi functional. In this set-up,  K\"ahler-Ricci solitons and the Kahler-Ricci flow play similar roles to the extremal metrics and the Calabi flow.

We now introduce a crucial notion in the theory, the 
 {\it Futaki invariant} \cite{kn:Fut}. Initially we define this for  a K\"ahler manifold $X$ with an $S^{1}$ action.

\ 

{\bf Differential geometric approach} Consider an $S^{1}$-invariant metric
$\omega$ in the given K\"ahler class. The $S^{1}$ action is generated by
a Hamiltonian function $H$. The Futaki invariant is defined to be
\begin{equation} {\rm Fut}=  \int_{X} H (S-\hat{S}) \frac{\omega^{n}}{n!}. \end{equation}

It does not depend on the choice of metric $\omega$. 
In particular, if the Futaki invariant is not zero there can be no CSCK metric.
Extremal
metrics are the appropriate generalisation of CSCK to the case of manifolds
with non-trivial automorphisms and non-zero Futaki invariant.

\

{\bf Algebro-geometric approach}
Assume that the K\"ahler class corresponds to an integral polarisation
so we have an ample line bundle $L\rightarrow X$ with $c_{1}(L)= 2\pi [\omega]$ and fix a lift of the $S^{1}$-action to $L$.

Let $d_{k}= {\rm dim} H^{0}(X,L^{k})$. For large $k$ it is given by the Hilbert
polynomial, of degree $n$. There is an induced $S^{1}$-action on $H^{0}(X,L^{k})$: let $w_{k}$ be the weight of the action on the highest exterior power.  For
large $k$ it is given by a Hilbert polynomial, of degree $n+1$ (as explained in \cite{kn:SKD2.5} for example).   

Set $F(k)= \frac{w_{k}}{k d_{k}}$ so for large $k$ we have an expansion

$$  F(k) = F_{0} + k^{-1} F_{1} + k^{-2} F_{2}+ \dots. $$
The Futaki invariant is the co-efficient $F_{1}$.

\

The equivalence of the differential geometric and algebro-geometric definitions comes from the equivariant
Riemann-Roch formula.
The  algebro-geometric viewpoint has the advantage that it extends immediately to the case when $X$
is a singular variety,  or even a scheme. 

\subsection{The YTD conjecture}

This is the analogue of the Kobayashi-Hitchin conjecture for the case of CSCK metrics and to formulate it we need to specify the notion of stability which is called {\it K-stability}.

Fix $(X,L)$ as above. We consider {\it equivariant degenerations} (or {\it
test configurations})
$   {\cal L}\rightarrow {\cal X}\rightarrow \bC $.
Here \begin{itemize}
\item $\pi:{\cal X}\rightarrow \bC$ is a flat family with $\pi^{-1}(1)=X$;
\item  the line bundle ${\cal L}\rightarrow {\cal X}$ is ample on the fibres
and the restriction of ${\cal L}$ to $\pi^{-1}(1)$ is isomorphic to some power $L^{m}$;
\item  there is a $\bC^{*}$-action on the whole set-up.
\end{itemize}

  We define the Futaki invariant of a degeneration ${\cal X}$ to be the invariant
we defined above for the  central fibre $\pi^{-1}(0)$, polarised by
the restriction of ${\cal L}$, and
we say that $(X,L)$ is {\it K-stable} if ${\rm Fut}({\cal X})>0$ for all
non-trivial test configurations ${\cal X}$ as above.
There are some technicalities in the precise definition of \lq\lq non-trivial'', as pointed out by Li and Xu \cite{kn:LX}.  One can require that the total space ${\cal X}$ be normal and that ${\cal X}$ is not the product $X\times \bC$.

Then the (so-called) YTD conjecture is that $(X,L)$ admits a CSCK metric if and
only if it is K-stable. 

\

{\bf Remarks}

\begin{enumerate} \item It would fit better with the terminology used in 1.2 above to call this \lq\lq K-polystability'', as is done by some authors. 
\item In the case of Fano manifolds and K\"ahler-Einstein
metrics a version of this conjecture was proposed by Yau around 1990.
\item The definition of K-stability (again in the  K\"ahler-Einstein situation),
for the case when the central fibre is smooth or mildly singular was given
by Tian in 1996 \cite{kn:Ti}.  
\item The statement of the conjecture can be extended to include the case
of extremal metrics \cite{kn:Sz0}. 
\item   The \lq\lq easy'' side of the conjecture is generally seen as the fact that CSCK implies K-stability and the results here are here are relatively complete.  Donaldson proved that CSCK implies K-semistability \cite{kn:SKD2.5} and this was refined to K-stability by Stoppa\cite{kn:St}, assuming that the automorphism group is finite (i.e. $H^{0}(TX)=0$). In the Fano case the complete result was proved by Berman (with an extension to the singular case) \cite{kn:B}. Similarly there are complete results on the uniqueness of CSCK metrics (modulo holomorphic automorphisms). For K\"ahler-Einstein metrics this was proved by Bando and Maubuchi \cite{kn:BM}. In the CSCK case, uniqueness would follow immediately from the convexity of the Mabuchi functional if one knew that any two metrics could be joined by a smooth geodesic. While this is known {\it not} to be exactly true, variants of the argument can be made to work. This was done by Chen \cite{kn:Chen} in the case when $c_{1}<0$, by Donaldson in the projective case using a finite dimensional approximation and assuming that $H^{0}(TX)=0$ \cite{kn:SKD1} and by Chen-Tian \cite{kn:CT}. The most complete and decisive results have been obtained by Berman and Berndtsson \cite{kn:BB}; see also \cite{kn:BDL}, \cite{kn:CS}.

\item The \lq\lq correct'' conjecture may be a little different (see 4.3 below).

\item There is a circle of ideas and results relating $K$-stability to the better-established notions in algebraic geometry of 
Chow stability and Hilbert stability.  These proceed via projective embeddings
$X\rightarrow \bP(H^{0}(L^{k})^{*})=\bP^{N_{k}}$ and applying Geomtric Invariant theory to the action of $SL(N_{k}+1, \bC)$ on the appropriate Chow variety and Hilbert scheme. 
 As $k\rightarrow \infty$
there are asymptotic relations, connected to the subject of  geometric quantisation,  between
these finite-dimensional pictures and the infinite dimensional picture. For example it can be proved in this way if $X$ has a CSCK metric and finite automorphism group then it is Chow stable for sufficiently large $k$ \cite{kn:SKD1}. But we will not go into this aspect in detail here. 
\end{enumerate}

\

 We will now attempt to motivate, informally, the definition of K-stability and hence the YTD conjecture. From one point of view, we know that a CSCK metric
corresponds to a critical point (in fact minimum) of the Mabuchi
functional) ${\cal F}$ on the space ${\cal H}$ of
K\"ahler metrics. Very roughly:
\begin{itemize}
\item  we expect that if there is no minimum then a minimising sequence will
tend to a \lq\lq point at infinity'' in ${\cal H}$;
\item there should be a numerical criterion which tells us which of the points
at infinity are \lq\lq destabilising''
i.e. whether the functional ${\cal F}$ decreases as we approach that point
at infinity.
\end{itemize}

In this vein the YTD conjecture can be thought of saying that the relevant points at infinity are \lq\lq algebro-geometric objects'', in fact the central fibres of test configurations, and that the sign of the Futaki invariant gives the appropriate numerical criterion. But we emphasise that this is just a motivating picture, which does not in itself make any progress towards a proof. \lq\lq Points at infinity'' in ${\cal H}$ have no {\it a priori} meaning and ${\cal H}$ is not even locally compact. 

From another point of view, we can compare with the finite dimensional situation where we have the Hilbert-Mumford criterion for stability in terms of 1-parameter subgroups $\lambda: \bC^{*}\rightarrow K^{c}$ and the resulting weight $w(\lambda,v)$.
This can be described as follows. We consider the action of the 1-parameter subgroup  on the point $[v]$ in the projective space $\bP(V)$. This has a well-defined limit
   $$   p = \lim_{t\rightarrow 0} \lambda(t) [v], $$
   which is a fixed point of the $\lambda$ action on $\bP(V)$. There is no loss of generality in assuming that $\lambda$ is the complexification of a circle subgroup of $K$, with generator $\xi\in {\rm Lie}(K)$.
   Then we have the formula
   \begin{equation} w(\lambda, v)= \langle \mu(p), \xi\rangle, \end{equation} where 
   $\xi \in {\rm Lie} (K^{c})$ is the generator of $\lambda$. In fact this is the weight of the action of $\lambda$ on the fibre of the tautological line bundle over $p$. 
Now we try to take this over to the case of the action of ${\cal G}$ on ${\cal J}_{{\rm int}}$, so we consider a circle subgroup of ${\cal G}$ generated by a Hamiltonian $H$. A fixed point in ${\cal J}_{{\rm int}}$ corresponds to a K\"ahler structure with a circle action,  and the formula (3) goes over to the formula (2) for the Futaki invariant in the smooth case. The YTD conjecture can then be viewed as positing
\begin{enumerate}
\item Test configurations can be thought as corresponding (at least roughly) to 1-parameter subgroups in ${\cal G}^{c}$.
\item We allow singular central fibres, so actually we move outside the space ${\cal J}_{{\rm int}}$, as we have defined it.
\item With these understandings the analogue of the Hilbert-Mumford criterion for the existence of the zero of the moment map is true.
\end{enumerate}

With regard to the first item, recall that in the finite-dimensional situation the geodesics in the symmetric space $K^{c}/K$ correspond to analytic 1-parameter subgroups in $K^{c}$ which are the complexification of $1$-parameter subgroups in $K$. So one can attempt to make sense of 1-parameter subgroups in ${\cal G}^{c}$ as  {\it geodesic rays} in the space of metrics ${\cal H}$. This geodesic equation makes perfectly good sense: it is a version of the homogeneous complex Monge-Amp\`ere  equation and various results in the direction of item (1) have been established, beginning with Phong and Sturm \cite{kn:PS}. 

\section{Proofs of cases of the YTD conjecture }

\subsection{Initial discussion}

The YTD conjecture, for general polarised manifolds $(X,L)$ remains a conjecture---and one which does not seem likely to be established in the near future. The constant scalar curvature equation is a 4th. order non linear PDE  and the existence theory is very limited at present.     There are results on small deformations by Sz\'ekelyhidi \cite{kn:Sz1.5} and others.  There are other results of an asymptotic nature  considering blow-ups and fibred manifolds. In the first case one considers a set of points in a CSCK manifold and  the existence problem for metrics on the blow-up, in a K\"ahler class where all the exceptional fibres are small \cite{kn:AP}. (Combining this  blow-up theory with other results, Shu showed that in each deformation class of K\"ahler surfaces there is   a manifold admitting an extremal metric \cite{kn:Shu}.)  In the second case, one considers a vector bundle $E\rightarrow X$ over a CSCK manifold and metrics on the projectivisation $\bP(E)$, in a K\"ahler class where the fibres are small. The CSCK equations are then related to the Hermitian Yang-Mills equation and the K-stability of $\bP(E)$ to the stability of $E$ \cite{kn:YH}\cite{kn:RossThomas}. There are also results for other  fibred manifolds \cite{kn:JF}.

The situations mentioned above are all of a perturbative nature and the proofs  are based on implicit function theorems. Beyond this there are two main cases where stability enters in an essential way and where the YTD conjecture has been established: toric surfaces and Fano manifolds. We will discuss these in the next two subsections. 

 \subsection{Toric manifolds}
 
 We consider a polarised toric manifold $(X,L)$ of complex dimension $n$. Thus the complex torus $T^{c}=(\bC^{*})^{n}$ acts on $(X,L)$ and there is a dense, free, open orbit in $X$. This data defines a convex polytope $P\subset \bR^{n}$, which is the convex hull of a finite set of integral points and which satisfied the \lq\lq Delzant condition''. Conversely, any such polytope yields a polarised toric manifold. This correspondence can be developed in many ways. In terms of K\"ahler geometry, we consider K\"ahler metrics on $X$ which are invariant under the real torus $T=(S^{1})^{n}$. Restricted to the open orbit, such a metric is given by a K\"ahler potential which can be viewed as a convex function on $T^{c}/T= \bR^{n}$. The Legendre transform of this function is then a convex function $u$ on $P$. This approach was developed by Guillemin \cite{kn:Guillemin} and Abreu \cite{kn:Ab} and corresponds exactly to the toric case of that described in 1.2 above. The manifold $X$ is regarded as a compactification of ${\rm int} P \times T$ with the fixed symplectic form $\sum dx^{a} d\theta_{a}$ where $x^{a}$ are standard co-ordinates on $\bR^{n}$ and $\theta_{a}$ are standard angular co-ordinates on $T$. Then the convex function $u$ defines a complex structure on $X$, specified by saying that the complex $1$-forms  $d\theta_{a} + i \sum_{b} u_{ab}d x^{b}$ have type $(1,0)$. Here $u_{ab}$ denotes the second derivative
  $ \frac{\partial^{2} u}{\partial x^{a} \partial x^{b}}$. The function $u$ is required to satisfy certain boundary conditions on $\partial P$---roughly speaking, $u$ should behave like $d\log d$ where $d$ is the distance to the boundary. These boundary conditions mean that the complex structure extends smoothly to the compact manifold $X$. While these complex structures---for different symplectic potentials $u$---are different they are all isomorphic. The relation between the two points of view, fixing either the complex structure or the symplectic structure, becomes the classical Legendre transform  for convex functions.

  Each codimension-$1$ face of $P$ lies in a hyperplane which contains an integer lattice. This lattice defines a Lebesgue measure on the hyperplane. Putting these together, we get a natural measure $d\sigma$ on $\partial P$. A crucial object in the theory is then the linear functional ${\cal L}$ on functions $f$ on $P$, 
$$\cL  f = \int_{\partial P} f  d\sigma - A \int_{P} f d\mu. $$
Here the real number $A$ is
   $   {\rm Vol}(\partial P, d\sigma)/{\rm Vol}(P,d\mu)$,
   so that $\cL$ vanishes on the constant functions. The restriction of $\cL$ to the linear functions defines the Futaki invariant of $(X,L)$ and we can only have a constant scalar curvature metric if this vanishes, which we now assume. That is, we assume that the centre of mass of $(\partial P, d\sigma)$ coincides with the centre of mass of $(P,d\mu)$.

The metric defined by a convex function $u$ is
  \begin{equation}    \sum u_{ab} dx^{a} dx^{b} + u^{ab} d\theta_{a} d\theta_{b}, \end{equation}

where $\left( u^{ab}\right)$ is the inverse of the Hessian matrix $\left(u_{ab}\right)$. The scalar curvature is
  \begin{equation} S = -\frac{1}{2} \sum_{ab} \frac{\partial^{2} u^{ab}}{\partial x^{a} \partial x^{b}}. \end{equation}
If we have a constant scalar curvature metric then the constant is $A/2$, so the problem is to solve the fourth order nonlinear PDE
  \begin{equation}  \sum_{ab} \frac{\partial^{2} u^{ab}}{\partial x^{a} \partial x^{b}} = -A \end{equation}
for a convex function $u$ satisfying the boundary conditions and with this given  $A$. This has a variational formulation. We set
      \begin{equation}   {\cal F}(u)= - \int_{P} \log \det \left(u_{ab}\right) + \cL(u), \end{equation}
which is just the Mabuchi functional in this context.
Then the PDE problem is to find a minimum of the functional $\cF$ over all convex functions on $P$.

We now turn to algebraic geometry. 
   Suppose that $f$ is a convex {\it rational piecewise-linear
function}, i.e.
$ f=\max\lambda_{i}$ where $\{\lambda_{i}\}$ are a finite collection of affine-linear functions with rational co-efficients. Thus $f$ defines a decomposition of $P$ into convex rational polytopes, on each of which $f$ is affine-linear. Define $Q\subset P\times \bR$ to be $\{(x,y):  y\geq f(x)\}$.
This is a  convex polytope which corresponds
to an $(n+1)$-dimensional toric variety ${\cal X}$. This yields a degeneration
of $X$ in which the central fibre is typically reducible: the components of the central
fibre correspond to the pieces in the decomposition of $P$. For example if we take $n=1$, $P=[-1,1]\subset \bR$ and $f(x)= \vert x\vert$ then we get the degeneration of $\bC\bP^{1}$, embedded as a smooth conic in $\bC\bP^{2}$, into a pair of lines. 

Now the Futaki invariant
of such a degeneration ${\cal X}$ is just $\cL(f)$. So we say that $(X,L)$ is {\it toric K-stable} if $\cL(f)\geq 0$ for all such rational piecewise-linear  functions $f$, with equality if and only $f$ is affine-linear. The  YTD conjecture in this case is the conjecture that this is the necessary and sufficient condition for solving the PDE problem (6). (In fact there is a significant technical difficulty here, in showing that for toric manifolds toric K-stability is equivalent to K-stability as previously defined---{\it a priori} the manifold could be de-stabilised by a degeneration which is not compatible with the toric structure. This was resolved recently---for toric actions---by work of Codogni and Stoppa  \cite{kn:CS})

Some part of the connection between the PDE problem and the stability condition can be seen as follows. Suppose that there is a PL convex function $f$ with $\cL(f)<0$. Then we  choose a smooth convex function $\tilde{f}$ with $\cL(\tilde{f})<0$ and, given some reference potential $u_{0}$, we consider the 1-parameter family
$u_{s}= u_{0}+ s \tilde{f}$. The fact that $\tilde{f}$ is convex means that for $u_{s}$ is a symplectic potential for all $s\geq 0$. The slow growth of the logarithmic term in (7) means that
\begin{equation}     \cF(u_{s})\sim s \cL(\tilde{f}) \end{equation}
as $s\rightarrow \infty$. In particular the functional $\cF$ is not bounded below and there can be no minimum. This fits in with the general heuristic picture we outlined in 1.3 above. The metric on the space of  torus-invariant  K\"ahler metrics is flat and the geodesics are exactly the linear paths like $u_{s}$.

We note that for toric Fano manifolds the existence problem for K\"ahler-Einstein, metrics and more generally K\"ahler-Ricci soliton metrics, was completely solved by Wang and Zhu \cite{kn:WZ}. It is also possible to verify that the appropriate stability conditions are satisfied in this situation.

\

The toric version of the YTD conjecture in the case of toric surfaces, $n=2$,  was confirmed in a series of papers culminating in \cite{kn:SKD3}. The strategy of proof was a continuity method. Suppose that we have  any polytope $P$ (not necessarily having the Delezant property, or satisfying any rationality condition) with a measure $\partial \sigma$ on the boundary. Then we can define the linear functional $\cL$ and set up the PDE problem for a function $u$.  We introduce the stability condition, that $\cL(f)\geq 0$ for {\it all} convex functions $f$, with equality if and only if $f$ is affine-linear. Then part of the work is to show that in the case of a pair $(P,d\sigma)$ arising from a toric surface this notion coincides with that defined before. That is, a destabilising function $f$ can be taken to be  piecewise linear and rational. With this alternative definition in place one can consider $1$-parameter families $(P_{t}, d\sigma_{t})$ and the main task is to show that if the stability condition holds for all $t$ in a closed interval $[0,1]$ and if there is a solution $u^{(t)}$ for $t<1$ then in fact there is a solution for $t=1$. This requires many steps. One set of arguments culminate in a uniform bound \begin{equation} \Vert u^{(t)} \Vert_{L^{\infty}} \leq C \end{equation}
for all $t<1$. From this we get a weak limit $u^{(1)}$ (initially from   a subsequence $t_{i}\rightarrow 1$) and the problem is to show that this is smooth, strictly convex and satisfies the boundary conditions.  This requires further estimates. For example a relatively simple estimate gives a uniform lower bound
\begin{equation}   \det\left( u^{(t)}_{ab}\right) \geq \epsilon>0. \end{equation}
One of the  special features of dimension 2 is that such a lower bound on the Monge-Amp\`ere function gives a strict convexity property of $u$---this fails in higher dimensions. In general, 
the estimates near the boundary are particularly difficult.

All of these results, for toric surfaces, have been extended to extremal metrics by B. Chen, A-M Li and L. Sheng \cite{kn:CLS}. The extremal equation in this setting is the same equation (6) but now with $A$ an affine-linear function on $\bR^{n}$.

  The underlying reason why the differential geometric theory
works   simply in the toric case is
in that the group of symplectomorphisms which commute with the
torus action is {\it abelian}. In the world of symplectic manifolds with  group actions there is a larger class of \lq\lq multiplicity free''  manifolds
having this property. These roughly correspond to \lq\lq reductive varieties'' in algebraic geometry. Much of the theory discussed above can be extended to this situation\cite{kn:AK}, \cite{kn:SKD5}. Recently Chen Han, Li, Lian and  Sheng have extended their analytical results to establish a version of the YTD conjecture for a class of these manifolds, in the  rank 2 case \cite{kn:CHLLS}.  A comprehensive treatment of the Fano case has been achieved  by Delcroix \cite{kn:Delcroix}.

\

To understand the existence proofs, it is important to understand what happens when the stability condition fails. Suppose then that the pair $(P_{1}, d\sigma_{1})$ corresponds to a K-semistable surface $X$, so that there is a non-trivial  convex function $f$ with $\cL(f)=0$. One can show that $f$ can be taken to have the form $\max (0,\lambda)$, where $\lambda$ is affine linear; so the resulting decomposition of $P_{1}$ has just two pieces $P=P^{+}\cup P^{-}$. This  corresponds to a degeneration where the central fibre has two components
$ Y^{+}\cup_{D}Y^{-}$. Now suppose that this pair is embedded in a 1-parameter family as above, where the stability condition holds for $t<1$. The uniform estimate(9) will then fail. The analytical result suggest that we can normalise $u^{(t)}$ in two different ways
$$   u_{+}^{(t)}= u^{(t)}+ c^{+}_{t} f  \ ,\  u_{-}^{(t)}= u^{(t)}+ c^{-}_{t} f, $$
 for $c^{\pm}_{t}\in \bR$, so that $u_{\pm}^{(t)}$ converge over $P^{\pm}$ to limits $u_{\pm}^{(1)}$,  and that these limits define complete CSCK metrics on the non-compact surfaces
$Y^{\pm}\setminus D$. Transverse to the divisor $D$, the picture should be modelled on the well-known degeneration of a family of hyperbolic Riemann surfaces into a pair of cusps. (We can suppose that the boundary between $P^{+},P^{-}$ is a line $x^{1}={\rm constant}$. Then the second derivative of $u^{t}_{11}$ becomes very large  and one sees from the formula (4) that the metric becomes large in the $x^{1}$ direction and small in the $\theta_{1}$ direction.)

\section{ K\"ahler-Einstein metrics on Fano manifolds}
\subsection{Initial discussion}

In this section we discuss  the proof of the Yau (or
YTD) conjecture in the case of  Fano manifolds  by Chen, Donaldson and Sun \cite{kn:CDS}. That is:

\

{\bf Theorem} {\em A Fano manifold $X$ admits a K\"ahler-Einstein metric if and
only if  $(X,K_{X}^{-1})$ is K-stable.}

\

We should emphasise again that there are many earlier results,  of Siu \cite{kn:Siu}, Nadel\cite{kn:Na}, Tian \cite{kn:Ti0} and others, proving that particular Fano
manifolds are K\"ahler-Einstein, using the theory of the $\alpha$-invariant
and log canonical threshold.

\

 The equation we want to solve is $\rho=\omega$ where
$\rho$
is the Ricci form. There are at least three strategies to attack this PDE problem.
\begin{enumerate} \item {\bf The continuity method}: for a parameter $s\in [0,1]$
and some fixed positive closed form $\rho_{0}$ try to solve
$$  \rho= s \omega + (1-s) \rho_{0},  $$
and hope to solve up to $s=1$ to get the K\"ahler-Einstein metric.
\item {\bf K\"ahler-Ricci flow}:  with any initial condition there is a solution
of the flow equation 
$$   \frac{\partial \omega}{\partial t} = \omega-\rho, $$
defined for all $t>0$. We hope to obtain a K\"ahler-Einstein metric as the limit when $t\rightarrow
\infty$. 
\item {\bf Cone singularities}: Fix a smooth divisor $D\in \vert -p K_{X}\vert$
for
some $p$. For $\beta\in (0,1]$ we try to find a K\"ahler-Einstein metric with
a cone  angle $2\pi \beta$ transverse to $D$ and hope to solve up to $\beta=1$
to get a smooth K\"ahler-Einstein  metric. 
\end{enumerate}

These have obvious similarities in spirit and also in the more technical aspects---the crucial thing is to obtain the appropriate limit which will of course involve using the stability condition (of course the same holds for the toric case discussed in the previous section).

 Approach (3) was the one taken by Chen, Donaldson, Sun. Subsequently other
proofs
were  given by Datar and Sz\'ekelyhidi  \cite{kn:Sz2}, \cite{kn:DS} using approach (1) and by  Chen and Wang \cite{kn:CW} and Chen, Sun and Wang \cite{kn:CSW} using approach (2). In the unstable case, Chen, Sun and Wang  show that  the  flow converges to a Ricci soliton metric on a (possibly singular) variety. They show that this variety is a degeneration of the central fibre of a destabilising test configuration. These other proofs have the advantage that they are compatible with group actions and lead to explicit new examples of manifolds admitting  K\"ahler-Einstein metrics. (For torus actions, this can now also be handled using the result of Codogni and Stoppa \cite{kn:CS} mentioned before.) We should also mention that there is now another proof (of a slightly different result) due to Berman, Boucksom, Jonsson \cite{kn:BBJ} which uses very different ideas.

In all three strategies (1), (2), (3) , the crucial thing is to be able to take a {\it limit} of
the metrics in our approximating scheme and to show that this limit is an
algebro-geometric object. The foundation for this is the existence of a deep convergence theory for
Riemannian metrics given suitable control of the Ricci tensor. This is what
makes the K\"ahler-Einstein problem more tractable than the general CSCK
one.  In the next two subsections we will outline some of the main ideas involved.

\subsection{Gromov-Hausdorff convergence}

Let $A,B$ be compact metric spaces. The Gromov-Hausdorff distance $d_{GH}(A,B)$
can be defined by saying that for $\delta>0$ we have $d_{GH}(A,B)\leq \delta$
if there is a metric on $A\sqcup B$ extending the given metrics on $A,B$
and such that $A$ and $B$ are $\delta-dense$. A fundamental theorem of Gromov runs as follows. 
Suppose $C,D>0$ are given and $(M_{i}, g_{i})$ is a sequence of Riemannian
manifolds of fixed dimension $m$ and with
\begin{itemize}\item ${\rm Ricci}(M_{i}, g_{i})\geq C$;
\item ${\rm Diam}(M_{i}, g_{i}) \leq D$.
\end{itemize}
Then there is a subsequence which converges, in the sense of $d_{GH}$, to
some limiting metric space.

The Ricci curvature enters the proof through volume comparison. For simplicity
suppose $C=0$. Bishop's Theorem states that for a point $p$ in a compact manifold $M$ with
${\rm Ricci}\geq 0$ the ratio
$$   \frac{{\rm Vol}(B_{p}(r))}{r^{m}} $$
is a decreasing function of $r$. This is a strong global form of the infinitesimal relation (1). The consequence here is that if the volume of $M$ is $V$ and the diameter
is $\leq D$ then 
\begin{equation}  {\rm Vol}(B_{p}(r)) \geq \frac{V}{D^{m}} r^{m} . \end{equation}

Now a packing argument shows that for any $\epsilon>0$ there is a fixed computable
number $N(\epsilon)$ (determined by $V$ and $D$),  such that $M$ can be covered by $N(\epsilon)$ balls of
radius $\epsilon$. Applying this to the $M_{i}$, the construction of the
subsequence and the Gromov-Hausdorff limit follows from an elementary argument.
(Take a sequence of approximations to the $M_{i}$ by {\it finite sets} and
use the fact that  any bounded sequence of real numbers has a convergent
subsequence along with a \lq\lq diagonal argument''.)

 While the Gromov-Hausdorff  limit is initially just a metric space a lot more
is known about it. Results of Cheeger and Colding \cite{kn:CG1}, Anderson \cite{kn:An} and Cheeger, Colding and Tian \cite{kn:CCT} establish that:
\begin{itemize}
\item if the Ricci tensors of the $M_{i}$ also have an upper bound and
the volumes of $M_{i}$ are bounded below then
 the limit $M_{\infty}$ is the union $R\sqcup S$ where the {\it regular
set} $R$ is a $C^{1,\alpha}$ Riemannian manifold 
and the {\it singular set} $S$ is closed, of Hausdorff codimension $\geq
4$;
\item at each point of $M_{\infty}$ there exist metric \lq\lq tangent
cones''.
\end{itemize}

\subsection{Gromov-Hausdorff limits, line bundles and algebraic geometry}

 To get to the main ideas we discuss a slightly different situation as in \cite{kn:DS}. Suppose that we have a polarised manifold $L\rightarrow X$ and a K\"ahler
class $2\pi c_{1}(L)$, a   Hermitian
metric on $L$ and the K\"ahler form defined by the  curvature of this metric.
A fundamental, but completely elementary , point is that 
taking a power $L^{k}$
of the line bundle corresponds to scaling the metric on $X$, so distances are scaled by $\sqrt{k}$.

Suppose that we have a sequence of such $(L_{i}, X_{i}, \omega_{i})$ with
a fixed diameter bound and bounded  Ricci curvatures.
By Gromov's Theorem, we can suppose there is a Gromov-Hausdorff  limit $Z$. The central problem is to relate this to algebraic geometry. A rough statement expressing this  is  that $Z$  is homeomorphic to a normal complex projective variety.
There are various more precise statements: for example (after perhaps taking a subsequence)
we can find a fixed $k$ and embeddings $\tau_{i}: X_{i}\rightarrow\bC\bP^{N}$
defined by the sections of $L_{i}^{k}$ such that the projective varieties
$\tau_{i}(X_{i})\subset \bC\bP^{N}$ converge to a normal variety $W$ which
is homeomorphic to $Z$.

The essential difficulty in proving this is to establish a \lq\lq partial $C^{0}$-estimate''.
Given a positive Hermitian line bundle $L\rightarrow X$ and a point $p\in
X$ we define $\rho(p,L)$  to be the norm of the evaluation map
$$    ev_{p}: H^{0}(X,L)\rightarrow L_{p}, $$
using the standard $L^{2}$ metric on the sections. 
The statement that $\rho(p,L)>0$ for all $p$ is the statement that the sections
of $L$ define  a map 
$$  \tau: X\rightarrow \bP(H^{0}(X,L)^{*}). $$

 A definite lower bound on $\rho(p,L)$ gives metric control of the map $\tau$.
The estimate in question is to show that there is some fixed $k_{0}$ and
a definite lower bound on $\rho(L_{i}^{k_{0}}, p)$ for all $p\in X_{i}$. This can be thought of as a quantitative form of the standard Kodaira theory for high powers of a positive line bundle.

The outline of the argument to achieve this is as follows. 

\begin{itemize} \item We work near a point $q$ of the Gromov-Hausdorff  limit $Z$. A very
small neighbourhood of $q$ is close to a neighbourhood in a tangent cone
$C(Y)$ (i.e the metric cone over a metric space $Y$).  
\item For a suitable open set $U$ in the regular part of $C(Y)$ we get 
approximately holomorphic and isometric embeddings $\chi_{i}:U\rightarrow
X_{i}$,  after rescaling and for large  $i$.
\item Over $U$, in the regular set of the cone, there is a canonical section
$\sigma$ of the trivial holomorphic line bundle. It has Gaussian decay:
$$  \vert \sigma\vert^{2} = e^{-r^{2}/2}, $$
where $r$ is the distance to the vertex of the cone. 
\item If we can trivialise (or \lq\lq approximately trivialise'') the bundle
$\chi_{i}^{*}(L_{i}^{k})$, for a suitable $k$, and if we have a suitable
cut-off function $\beta$ supported in $U$, we can transport $\beta \sigma$
to a $C^{\infty}$ \lq\lq approximately holomorphic'' section $\tilde{\sigma}$
of $L_{i}^{k}$. 
\item There is a well-developed \lq\lq H\"ormander technique''
for studying the projection of  $\tilde{\sigma}$ to the  holomorphic
sections.
\end{itemize}
 
The main work goes into the construction of suitable cut-off functions and
in analysing the holonomy of the line bundles. 
The upshot is that one constructs sections of $L_{i}^{k}$ which are sufficiently
controlled to give the required estimate. 
\subsection{Outline of the main proof (the YTD conjecture for Fano manifolds)}
Following the approach (3), using K\"ahler-Einstein metrics with cone singularities, the
programme of work is.
\begin{enumerate}\item Show that there is a solution for small cone angle.
\item Show that given one solution the cone angle can be slightly deformed.
\item Extend the discussion above to metrics with cone singularities to show
that if there is a sequence $\beta_{i}\rightarrow \beta_{\infty}>0$ with
solutions for each $i$ then we can take a Gromov-Hausdorff limit which is
naturally a normal projective variety $W$ and is the limit of projective
embeddings  of $X$.
\item If $W=X$ show that there is a solution for the limiting cone angle
$\beta_{\infty}$.
\item If $W\neq X$, construct a non-trivial test configuration ${\cal X}$
with central fibre $W$.
\item Show that ${\rm Fut}({\cal X})\leq 0$.
\end{enumerate}

In item (5) a key step is to establish that the automorphism group of $W$
is reductive. This is an extension of the standard Matsushima theorem to the singular case and depends upon sophisticated results from pluripotential theory \cite{kn:Bern}.
 Item (6) follows from  an extension of the definition of the Futaki invariant
to pairs.

An interesting feature is that this proof shows that to test K-stability for Fano manifolds it suffices to consider test configurations with normal central fibre. In this proof the fact emerges from the differential geometry and Riemannian convergence theory: it is bound up with the non-collapsing inequality (11). On the algebro-geometric side this fact was proved (a little earlier) by Li and Xu \cite{kn:LX}, using results from the minimal model programme.

\section{ Concluding discussion}
\subsection{Examples}

A major difficulty in this area is that it usually very hard to check whether a manifold is K-stable. This is one of the advantages of the toric case, where the criterion is relatively explicit. The situation for  Fano manifolds is much less satisfactory and there are few interesting examples known. Even in dimension 3, where the Fano manifolds are completely classified, it is not known exactly which of them are K-stable. This is an outstanding algebro-geometric problem and one can certainly hope that it will be much better understood by the time of the next  AMS summer algebraic geometry  meeting. Among recent developments we note the work of Ilten and S\"uss \cite{kn:IltenSuss} and Delcroix \cite{kn:Delcroix2}on manifolds with large symmetry groups and of Fujita and Odaka  \cite{kn:FOS} on a new \lq\lq $\delta$-invariant'' criterion. 

There is one interesting class of examples in dimension 3, which go back in this context to Tian \cite{kn:Ti}. These occur in the family of Fano 3-folds of type $V_{22}$ which can all be embedded in the Grassmannian of $3$-planes in $\bC^{7}$.  More precisely, let $U$ be a $7$-dimensional complex vector space and  $\Pi$ be a $3$-dimensional subspace of $\Lambda^{2} U^{*}$. Then we define a variety
$X_{\Pi}$ in ${\rm Gr}_{3}(U)$ to be set of  3-planes $P\subset U$ such that
$\omega\vert_{P}=0$ for all $\omega \in \Pi$. For generic $\Pi$ this is a smooth Fano 3-fold, so we have a family of manifolds parametrised by an open subset $\Omega\subset {\rm Gr}_{3}(U^{*})$. The problem is to identify the set $\Delta\subset \Omega$ of points $\Pi\in \Omega$ which define K-unstable manifolds $X_{\Pi}$. 

This is an interesting case  because it is known that neither $\Delta$ nor $\Omega\setminus \Delta$ are empty. We will recall some of the discussion of this from \cite{kn:SKD5}. Take $U$ to be the irreducible $7$-dimension representation $s^{6}$ of $SL(2,\bC)$. This is an orthogonal representation and the image of the Lie algebra under the action defines a 3-dimensional subspace $\Pi_{0}$ of $\Lambda^{2}U^{*}$. The corresponding variety $X_{0}=X_{\Pi_{0}}$ ---the Mukai-Umemura manifold---is smooth and clearly supports  an induced $SL(2,\bC)$-action. It can be shown that $X_{0}$ admits a K\"ahler-Einstein metric and hence is K-stable. 

One can use deformation theory to describe explicitly the intersection of $\Delta$ with a neighbourhood of the point $[\Pi_{0}] \in \Omega$.  
 The group $SL(2,\bC)$ acts on the deformation space $T=H^{1}(TX_{0})$ and as 
 a representation  of $SL(2,\bC)$ we have $T= s^{8}(\bC^{2})$. For each
small $\alpha\in
T$ we have a deformation $X_{\alpha}$. 
Then for  {\it small} $\alpha\neq 0$
the manifold $X_{\alpha}$ has a KE metric if and only if $\alpha$ is polystable (in
the GIT sense that we defined in Section 1.2) for the $SL(2,\bC)$ action. That is, if $\alpha$ is regarded as a section of a line bundle of degree $8$ over $\bC\bP^{1}$ then either
\begin{enumerate}
\item There is no zero of multiplicity $\geq 4$ (the stable case), or, 
\item There are two distinct zeros, each of multiplicity $4$.
\end{enumerate}
In the second case, there is a $\bC^{*}$-action on the  manifold $X_{\alpha}$. 
 There is a unique $SL(2,\bC)$-invariant divisor $D$ in the anticanonical linear system  $\vert -K_{X_{0}}\vert$. This is the variety of tangents to the rational normal curve in $\bP^{12}= \bP (s^{12})$ and its normalisation is $\bP^{1}\times \bP^{1}$.  The divisor $D$ is homeomorphic to $\bP^{1}\times \bP^{1}$ but has a cusp singularity along the diagonal. In suitable local co-ordinates $(x,y,t)$ on $X_{0}$ the divisor $D$ is defined by the equation
$y^{2}= x^{3}$. The local versal deformation of this cusp singularity is given by 
\begin{equation} y^{2}= x^{3} + a(t) x + b(t). \end{equation}
Globally, $a, b$ are sections of ${\cal O}(8), {\cal O}(12)$ respectively over the diagonal $\bP^{1}$. This has the following interpretation, at least at the infinitesimal level. As a representation of $SL(2,\bC)$ the deformation space  of the pair $(X_{0}, D_{0})$ is $s^{8}\oplus s^{12}$ with the $s^{8}$ factor corresponding to the deformations of $X_{0}$ (the space $T$ above) and $s^{12}$ to the deformations of $D_{0}$ within $X_{0}$. Then in the description by $(a,b)$ the term $a$ lies in  $s^{8}$ and $b$ in $s^{12}$. If $a$ is not stable, with a single zero of multiplicity $4$ at the point $t=0$ say, we can choose $b$ to have a zero of multiplicity greater than $6$ at $t=0$. We get a pair $(X_{\alpha}, D_{\alpha})$ whose orbit under an action of a 1-parameter subgroup $\bC^{*}\subset SL(2,\bC)$ converges to $(X_{0}, D_{0})$.  This leads to a plausible conjectural description of the unstable set $\Delta$.  Consider the surface singularity at $0$ in $\bC^{3}$ defined by the equation $f=0$ where  $$ f(x,y,t)= y^{2}- x^{3}- t^{4} x. $$
The versal deformation space is the quotient of $\bC[x,y,t]$ by the ideal generated by $f$ and the partial derivatives $f_{x}, f_{y}, f_{t}$,  and this has dimension $10$. So we expect that in the $22$ dimensional family of $K3$ surfaces arising from deformations of $D_{0}$ there should be a $12$ dimensional sub-family with singularities modelled on $\{f=0\}$. The deformation discussion suggest that these arise as surfaces which lie in unstable 3-folds $X_{\alpha}$, with a $6$-dimensional linear system of these singular surfaces for each such $X_{\alpha}$. Define   $\Delta'$ to be the locus of points $\Pi$ such that there is an anticanonical divisor in $X_{\Pi}$ which has  singularity of the form $f=0$. Then it seems plausible that 
  $X_{\Pi}$ is not K-stable if and only if it does not admit a $\bC^{*}$ action and if $\Pi$ lies in the closure of $\Delta'$. (We emphasise that the deformation discussion does not establish anything about the stability for points $\Pi$ far away from $\Pi_{0}$---{\it a priori} there could be quite different kinds of destabilising test configurations.)

\subsection{Connections with moduli theory}
For manifolds with $c_{1}< 0$ the existence problem for KE metrics is
 completely understood
following Aubin and Yau. But for canonically polarised {\it singular} varieties
 Odaka showed that K-stability is equivalent to the stability in the sense
of Koll\'ar, Shepherd-Barron and Alexeev \cite{kn:Od1} , which is formulated in terms of having at worst \lq\lq semi log canonical'' singularities.
There is a good generalisation of the notion of a K\"ahler-Einstein metric to the case of singular varieties (roughly, a metric on the smooth part defined by a bounded potential) \cite{kn:BBEGZ}. 
 Berman and Guenancia show that the existence of such a metric is equivalent to K-stability \cite{kn:BG}. That is,  the \lq\lq YTD conjecture'' is non-trivial and true for canonically
polarised varieties.

Now we turn to the Fano case. One can form  compactifications of  moduli spaces of
K-stable Fano manifolds (as topological
spaces)  using the the  Gromov-Hausorff topology defined via the K\"ahler-Einstein metrics.
There have been important  recent developments relating these to algebro-geometric
approaches. \begin{itemize}
\item Odaka, Spotti and Sun studied the case of surfaces \cite{kn:OSS} \item Spotti, Sun
and Yao extend the existence theory to K-stable ${\bf Q}$-Fano varieties. \cite{kn:SSY}.
\item Odaka constructed a compactified moduli space of K-stable ${\bf Q}$-Fano
varieties as an algebraic space \cite{kn:Od2}\item  Li, Wang and Xu show that the moduli
space of K-stable Fano manifolds is quasi-projective (and also construct the compactification as an algebraic space) \cite{kn:LWX} \end{itemize}
At the time of writing it seems not to be known if the Gromov-Hausorff compactification is projective, but one can hope for further developments in this direction.

\subsection{The CSCK and extremal cases}

Examples suggest that the definition above of K-stability may not be the
correct one for general polarised manifolds. Sz\'ekelyhidi has given a modified
definition in \cite{kn:Sz1}.  Given a polarised variety $(X,L)$ we write $R_{k}=H^{0}(X,L^{k})$ and form the graded ring $R= \bigoplus R_{k}$. Then we consider filtrations
    $$   \bC=F_{0}\subset F_{1} \subset F_{2} \dots \subset R, $$
    such that 
    \begin{itemize}
    \item $F_{i}. F_{j} \subset F_{i+j}$;
    \item $F_{i}= \bigoplus_{k} F_{i}\cap R_{k}$
    \item $R= \bigcup_{i} F_{i}$.
    \end{itemize}
    The {\it Rees algebra} of such a filtration ${\cal F}$ is the subalgebra of $R[t]$ given by:
    $$  {\rm Rees}({\cal F})= \bigoplus_{i} (F_{i} R)t^{i} . $$
 Sz\'ekelyhidi observes   if the Rees algebra is finitely generated the scheme
    $${\rm Proj}_{\bC[t]} {\rm Rees}({\cal F})$$ is a test configuration for $(X,L)$. Conversely a construction of Witt-Nystr\"om defines a filtration from any test configuration. So it is the same to talk about test configurations as filtrations with finitely generated Rees algebra. Sz\'ekelyhidi then extends the definition  of the Futaki invariant to general filtrations---not necessarily with finitely generated Rees algebra. In turn this leads to a new, and more restrictive,  notion of K-stability. For example in the toric case the vector space $R_{k}$ has  a standard basis labelled by the lattice points $\bZ^{n}\cap k P$. Given any convex function $f$ on $P$ we define  convex subsets
$$   P_{i,k}= \{ x\in P: f(x)\leq i/k\}. $$ 
and we define a filtration by taking $F_{i}\cap R_{k}$ to be the subspace corresponding to the lattice points in $k P_{i,k}$.  The Futaki invariant of this filtration  is given by the integral expression $\cL(f)$ we discussed before. 

It seems likely that this definition of Sz\'ekelyhidi gives the correct formulation of the YTD conjecture. For example in the toric case in dimension $n>2$ it seems likely that one has to consider general convex functions, not just rational piecewise linear ones. Perhaps there is a sensible way to define   \lq\lq points at infinity'' in the space of K\"ahler metrics and a precise relationship between  these  and filtrations.

\

\

We conclude with some remarks about the more differential-geometric and PDE aspects. Extending the existence theory to general CSCK and extremal metrics seems very hard at present because there  is no  analogue of the Riemannian convergence theory based
on control of just the scalar curvature. In the Fano case one is essentially concerned with  metrics having a positive lower bound on the Ricci curvature and this means that \lq\lq collapsing'' cannot occur---the Bishop inequality gives a lower bound (11) on the volume of any metric ball in terms of its radius. This is connected with the algebro-geometric fact that it suffices to test stability by degenerations to normal varieties  In general we have to consider such collapsing phenomena, as we see already in the toric case. (For example in a sequence of hyperbolic Riemann surfaces developing a long neck and degenerating into a pair of cusps,  the metric discs in the middle of the neck will have very small area.) From the algebro-geometric point of view, it is essential to consider degenerations to non-normal varieties. It seems reasonable to hope for progress in the case of higher dimensional toric manifolds or for general complex surfaces but there are many difficulties and  important problems for the future.



\end{document}